\documentclass[12pt]{article}
\input amssym.def
\input amssym.tex

 \DeclareSymbolFont{AMSb}{U}{msb}{m}{n}
\DeclareSymbolFontAlphabet{\mathbb}{AMSb}
\renewcommand{\Bbb}[1]{\mathbb #1}
\DeclareSymbolFont{Frak}{U}{euf}{m}{n}
\DeclareSymbolFontAlphabet{\mathfr}{Frak}
\renewcommand{\frak}[1]{\mathfr #1}

\begin{document}

\def\hook{\hbox to 12.5pt{\vbox{\vskip 6pt\hrule width 7pt height .5pt}
         \kern -4pt\vrule height 7pt width .5pt\hfil}}
\def\zbar{{\bar z}}
\def\wbar{{\bar w}}
\def\C{\Bbb C}
\def\R{\Bbb R}
\def\H{\Bbb H}
\def\P{\Bbb P}

\def\I{\bf I}
\def\J{\bf J}
\def\K{\bf K}

\def\i{\hat i}
\def\j{\hat j}
\def\k{\hat k}

\def\hq{{\hat q}}
\def\hpi{{\hat \pi}}
\def\hI{{\hat I}}
\def\cQ{{\cal Q}}
\def\cP{{\cal P}}

\def\pd#1#2{{\partial#1\over\partial#2}}
\def\pdb#1{{\partial \over \partial#1}}

\newtheorem{name}{Definition}
\newtheorem{thm}{Proposition}
\title{Generalized Symplectic Manifolds}
\author{Daniel Cartin  \\
\small
{\it Department of Physics} \\
\small
{\it Pennsylvania State University} \\
\small
{\it University Park, PA, 16802}}
\date{24 October 1997}

\maketitle

\begin{abstract}
This paper uses a generalization of symplectic geometry,
known as $n$-symplectic geometry and developed by
Norris\footnote{See~\cite{norris1}, as well as similar work in~\cite{awane,
deleon1, deleon2}.}, to find observables on three-dimensional manifolds. It
will be seen that
for the cases considered
, the $n$-symplectic observables are derivable from the symplectic observables
of $\C^2$. The quantization of these observables, as well as those on the
frame
bundle of $\R^n$, is also examined.

\end{abstract}

\section{Introduction}

One of the more vexing problems in mathematical physics in the process of
quantization. If one has a quantum theory, then by defining a kind of limiting
process, one can obtain a corresponding classical theory. What is not clear,
however, is if, by adding additional information to a classical theory, one
can find a consistent quantization of the system. Many methods have been used
to explore this question of a proper quantization procedure. One is that of
geometric
quantization, using a prescription of particle dynamics based on symplectic
geometry, which encodes Hamilton's equations of motion on the phase space.
The natural idea,
first suggested by Dirac~\cite{dirac}, is that the observables could be
quantized by mapping the anti-symmetric Poisson bracket to a commutation
relation of their operators. However, the failure of this procedure
is notorious, as witnessed by the Groenewold-van Hove theorems on 
$T^*S^1$~\cite{gotay:cyl}, $S^2$~\cite{gotay:sph}, and $\R^{2n}$
\cite{groen, vh1, vh2}.

So the question is how to
generalize this procedure. One might wonder
at the use of an inherently Newtonian theory, since one would also like to
include such ideas as the observer and the reference frame into a relativistic
model. This is the idea behind a generalization of symplectic geometry, known
as $n$-symplectic geometry. For symplectic geometry, the prototype space is the
cotangent bundle
$T^*M$, which can be thought of as position and momenta in a {\it given}
reference frame. In $n$-symplectic geometry, the bundle
of linear frames $LM$ is the basis of study. This
principal $GL(n,\R)$ bundle takes into account the lack of a preferential
reference frame, and instead looks at all of them at once. The general linear
group then represents the symmetries of particles in a physical system. In
addition, it can be seen that, since the cotangent bundle is an associated
bundle to $LM$, there is a map of the observables and Hamiltonian vector
fields on $LM$ to $T^*M$, so $n$-symplectic geometry is a covering theory for
the Hamiltonian theories of both particles and fields~\cite{fln}. Yet
not all systems will possess the full $GL(n,\R)$ symmetry. An example of
this is forming an orthogonal frame bundle, where the frames transform under
$SO(p,q)$, and represent the presence of an orientation and a
metric to a particular
manifold. To include this into the $n$-symplectic geometry, we can look not
only at linear frames, but any kind of principal fiber bundle, with some
symmetry
group $G$. 

This paper examines the more general cases of $n$-symplectic
manifolds. We study examples of low-dimensional manifolds that are not
symplectic, but are principal bundles. In Section 2, we go through
the basic results of $n$-symplectic geometry on the frame bundle. The two
basic
differences between symplectic and $n$-symplectic geometry are the
$\R^n$-valued
$n$-symplectic form and the group action on the $n$-symplectic manifold.
Section 3 looks at the quantization of the observables obtained via
$n$-symplectic geometry on the frame bundle of $L(\R^n)$. Although these
observables have the same functional dimension as those of symplectic
geometry,
we can find an algebra whose quantization map does not seem
to fall prey to the same kind of inconsistencies as the symplectic case.
In Section 4, the 2-symplectic geometry of the trivial bundle $\cP \simeq
\R^2 \times S^1$ is studied, and we obtain the vector-valued
observables. Some of these functions on $\cP$ can be seen to arise from the
symplectic geometry of $\C^2$, which, when subject to a
constraint, gives the observables of the 3-sphere, the subject of Section 5.
This restriction to $S^3$ gives observables which are functions of the three
spin variables $x^i$ which satisty the bracket relation $\{x^i, x^j \} =
\epsilon^{ijk} x^k$. Finally, some conclusions and directions for further work
are given in Section 6.

\section{Review of $n$-symplectic geometry}

\subsection{Motivation}

The canonical example of the $n$-symplectic manifold is that of the
frame bundle\footnote{In~\cite{deleon2}, the $n$-symplectic structure is
induced from that of the 1-jet bundle (see Example
3.2). We take the frame bundle as the canonical example because it is
has a physical motivation as a model for particle dynamics. See
also~\cite{fln}.},
so the question is whether this formalism can be generalized to other principal
bundles, and distinguished from the quantization arising from symplectic
geometry. This section gives a brief survey of $n$-symplectic geometry, the
details of which 
can be found elsewhere~\cite{norris1, norris2}. This description will be of
$n$-symplectic geometry on the prototype manifold, the bundle of linear frames,
a good place to explain as well as motivate the formalism. Later in this paper,
we will begin applying the ideas from $LM$ onto other principal fiber bundles.

We start with an $n$-dimensional manifold $M$, and let $\pi: LM \to M$ be the
space of linear frames over a base manifold $M$, the set of pairs
$(m, e_k)$, where $m \in M$ and $\{e_k\}, k=1, \cdots, n$ is a linear frame at
$m$.  This gives $LM$ dimension $n(n+1)$, with $GL(n,\R)$ as the structure
group acting freely on the right. We define local coordinates on $LM$ in terms
of those on the manifold $M$ -- for a chart on $M$ with coordinates $\{x^i\}$, 
let
\[
q^i (m, e_k) = x^i \circ \pi (m, e_k) = x^i(m) \qquad
\pi^i _j (m, e_k) = e^i \left( {\partial \over \partial x^j} \right)
\]
where $\{e^j\}$ denotes the coframe dual to $\{e_j\}$. These coordinates are
analogous to those on the cotangent bundle, except, instead of a single
momentum coordinate, we now have a momentum frame.  We want to place some kind 
of structure on $LM$, which is the prototype of $n$-symplectic geometry that is
similar to symplectic geometry of the cotangent bundle $T^*M$. The structure
equation for symplectic geometry
\[
df = - X \hook d\vartheta
\]
gives Hamilton's equations for the phase space of a particle, where
$\vartheta$ is the canonical symplectic 2-form. There is an
naturally defined $\R^n$-valued 1-form on
$LM$, the soldering form, given by
\[
\theta (X) \equiv u^{-1} [\pi_* (X)] \quad \forall X \in T_u LM
\]
where the point $u = (m, e_k) \in LM$ gives the isomorphism $u: \R^n \to
T_{\pi(u)}M$ by $\xi^i r_i \to \xi^i e_i$, where $\{r_i\}$ is the standard
basis of $\R^n$. The $\R^n$-valued 2-form $d\theta$
can be shown to be non-degenerate, that is,
\[
X \hook d\theta = 0 \Leftrightarrow X = 0
\]
where we mean that each component of $X \hook d\theta$ is identically zero.
Finally, since there is also
a structure group on $LM$, there are also group transformation properties. Let
$\rho$ be the standard representation of $GL(n,\R)$ on $\R^n$. Then it can
be shown that the pullback of $d\theta$ under right translation by $g \in
GL(n,\R)$ is $R^*_g \ d\theta = \rho(g^{-1}) \cdot d\theta$.

Thus, we have an $\R^n$-valued generalization of symplectic geometry, which
motivates the following definition.
\begin{name}
Let $P$ be a principal fiber bundle with structure group
$G$ over an m-dimensional manifold $M$. Let $\rho: G \to GL(n,\R)$
be a linear representation of $G$. An {\bf n-symplectic structure}
on $P$ 
is a $\R^n$-valued 2-form $\omega$ on $P$ that is (i) closed and
non-degenerate, in the sense that
\[
X \hook \omega = 0 \Leftrightarrow X = 0
\]
for a vector field $X$ on $P$, and (ii) $\omega$ is equivariant, such that
under the right action of $G$, $R^* _g \ \omega = \rho(g^{-1}) \cdot \omega$.
The  pair $(P,\omega)$ is called an {\bf n-symplectic manifold}.
\end{name}
Here, we have modeled $n$-symplectic geometry after the frame bundle by
defining
the general $n$-symplectic manifold as a principal bundle. There is no reason,
however, to limit ourselves to this, since we can let $P$ be any manifold with
a group action defined on it. One example of this would be to look at the
action of the conformal group on $\R^4$. Since this group is locally isomorphic
to $O(2,4)$, which is not a subgroup of $GL(4,\R)$, then forming a
$O(2,4)$ bundle over $\R^4$ cannot be thought of as simply a reduction of
the frame bundle. Another example is one that appears later in this paper,
namely, $\C^2$ with a $U(1)$ group action. This manifold has no principal
bundle structure, although $\C^2 - \{$origin$\}$ is a principal $U(1)$ bundle.
It is shown that this 2-symplectic geometry is actually just symplectic
geometry in disguise. Finally, note that in the definition, the dimension $n$
of the
representation space of $G$ is not necessarily equal to the dimension $m$ of
the base manifold $M$, although for the frame bundle and the examples
considered in this paper, $n = m$.

For $n$-symplectic manifolds other than $LM$, we will want to have the same
structure
equation, yet there will not be an obvious closed 2-form, such as the exterior
derivative of the soldering
form. Note that on the frame bundle, the $n$-symplectic form in local coordinates
is
\begin{equation}\label{struct}
d\theta = d\theta^i r_i = (d\pi^i _j \wedge dq^j) r_i
\end{equation}
Each component $d\theta^i$ of this 2-form is non-degenerate on a
$2n$-dimensional submanifold of $LM$, and so is a {\it symplectic} 2-form on
that submanifold. This paper will use this fact to construct 2-symplectic forms
using contact forms on three-dimensional manifolds.

\subsection{Observables on $LM$}

In local coordinates on $(LM, d\theta)$, we have the $n$-symplectic structure
equation given by
\begin{equation}\label{struct1}
df^{i_1 i_2 \cdots i_p} = -p! X^{i_1 i_2 \cdots i_{p-1}} _f \hook d\theta^{i_p}
\end{equation}
with the components $f^{i_1 i_2 \cdots i_p}$ of an $\otimes ^p \R^n$-valued
function $f$ and its associated Hamiltonian vector field, the $\otimes^{p-1}
\R^n$-valued $X_f$. Unlike symplectic geometry, the soldering form
transforms tensorially under right translations $R_g$ of the group $GL(n,\R)$, 
where $R^*_g \theta = g^{-1} \cdot \theta$ for $g \in GL(n,\R)$. Because of
this,
not every $\R^n$-valued function on $LM$ is compatible with the above
structure equation, as opposed to the fact that all smooth $\R$-valued
functions are allowable symplectic observables. First, we shall look at the
solutions to the equation
\[
df^i = - X \hook d\theta^i
\]
If we write our Hamiltonian vector field $X$ as
\[
X = X^a {\partial \over \partial q^a} + X^a _b {\partial \over \partial \pi^a
_b}
\]
then our two equations for the components of $X$ are similar to those of
symplectic geometry, namely,
\begin{equation}\label{t-one}
{\partial f^i \over \partial q^a} = - X^i _a \qquad
{\partial f^i \over \partial \pi^b _c} = X^c \delta^i _b
\end{equation}
It is this second equation for $X^c$ that gives the difference of symplectic
and $n$-symplectic observables. If we take the derivative of this equation by
$\pi^r _s$, then
\[
{\partial^2 f^i \over \partial \pi^b _c \partial \pi^r _s} =
{\partial X^c \over \partial \pi^r _s} \delta^i _b =
{\partial X^s \over \partial \pi^b _c} \delta^i _r
\]
If we then sum over the indices $i$ and $r$, then we find that, for $n >1$,
\[
{\partial X^s \over \partial \pi^b _c} = 0
\]
and the observables solving  $(\ref{t-one})$, designated $T^1$, are of the
form
\[
f^i = F^j (q) \pi^i _j + G^i(q)
\]

Now we look for solutions of  $(\ref{struct1})$ for $p>1$. The key point is
that, as this structure equation stands, the solutions would be the same as
functions in $T^1$, except that they would be $\otimes^p \R^n$-valued instead
of just $\R^n$-valued. For instance, the observables for $p=2$ would look like
\[
f^{ij} = {\tilde F}^{ib} (q) \pi^j _b + {\tilde G}^{ij} (q)
\]
There needs to be some kind of symmetry condition on the
structure equation. Solving the equation
\[
df^{i_1 i_2 \cdots i_p} = -p! X^{(i_1 i_2 \cdots i_{p-1}} _f \hook
d\theta^{i_p)}
\]
where the parentheses denotes symmetrization over the indices and $f$ is a
function totally symmetric on its indices, gives a space
$SHF^p$ of functions which, in local coordinates, are $\otimes^p _s
\R^n$-valued degree $p$ polynomial in the generalized momenta $\pi^i _j$
on $LM$\footnote{Note that the Hamiltonian vector field of observables in
$SHF^p$ is given by this
equation modulo a vector field $Y$ solving the kernel equation
\[
0 =  -p! Y^{(i_1 i_2 \cdots i_{p-1}} \hook d\theta^{i_p)}
\]
which gives us an {\it equivalence class} of Hamiltonian vector fields. It can
be shown~\cite{norris1} that if one always symmetrizes (or anti-symmetrizes,
for $AHF^p$) in the definition of quantities such as the Poisson bracket, then
the definition is independent of the choice of representative.}.
An example for $p=2$ is the observable
\[
g^{ij} = A^{(ab)} \pi^i _a \pi^j _b + 2 B^{c(i} \pi^{j)} _c + C^{(ij)}
\]
where $g^{ij} = g^{ji}$.
There is a similar observable if we anti-symmetrize over all indices, giving
us a space denoted $AHF^p$, but since the discussion is similar to the
symmetric case, we shall only deal with the latter. Also, notice that the
first term of $g^{ij}$ is equivariant under the action of $GL(n,\R)$. This
gives us a special space $ST^p$ of all homogeneous degree $p$ polynomials.
These observables are associated to symmetric tensors on the base manifold.
For example, suppose we have a vector field $\vec f$ on $M$.
We have the isomorphism $u^{-1} : T_{\pi(u)} \to \R^n$ for each point $u = (m,
e_k)$ in $LM$, which is given by
\[
u^{-1} (v) = e^i (v) r_i
\]
for a vector $v \in T_{\pi(u)}$. So, we define the function $f$ on $LM$ by
$f(u) = u^{-1}({\vec f} (\pi(u)))$. This gives the relation
\[
{\vec f} = A^i \pdb{x^i}\ \Leftrightarrow\ f = A^i \pi^j _i r_j
\]
This process can be suitably generalized for all symmetric tensors, by
extending $u^{-1}$ to a map from tensor products of $T_{\pi(u)}$ to tensor
products of $\R^n$.

One might question symmetrizing over all the indices, since, for $p>2$, the
general observable is not simply the sum of the symmetrized and
anti-symmetrized functions. For instance, we might look at an observable of
the form $f^{ijk} = f^{i(jk)}$. However, solving the equation
\[
df^{ijk} = -3! X^{i(j} _f \hook d\theta^{k)}
\]
would give us an observable which is quadratic, not cubic, in
the generalized momenta, while solving
\[
df^{ijk} = -3! X^{(jk)} _f \hook d\theta^i
\]
would give us an observable affine linear in $\pi^a _b$. This occurs for the
same reason going to the p=2 case without any kind of symmetry on the indices
does -- the structure equation does not take into account those indices which
are not included in the symmetry of the function $f$.

To form a Poisson algebra, we take the direct sum of all the $SHF^p$, called
$SHF$, and give it a Poisson bracket
\[
\{f,g\} = p! X^{(i_1 i_2 \cdots i_{p-1}} _f (g^{j_1 j_2 \cdots j_q)})
\]
where the symmetrization is necessary to get another element of $SHF$. Thus,
the Poisson bracket of elements in $SHF^p$ and $SHF^q$ give an element of
$SHF^{p+q-1}$. This bracket has the same properties as its symplectic
counterpart, namely,
\begin{eqnarray*}
&(1)& \{f \otimes_s g, h\} = f \otimes_s \{g, h\} + \{f, h\} \otimes_s g \\
&(2)& \{f, g\} = - \{g,f\} \\
&(3)& \{f, \{g,h\} \} +  \{g, \{h,f\} \} +  \{h, \{f,g\} \} = 0
\end{eqnarray*}
Despite the fact we are symmetrizing over the indices, the bracket is still
anti-symmetric, because of the anti-symmetry of the $n$-symplectic form. If we
look solely at the tensorial observables $ST$, then this bracket is the
version on $LM$ of the Schouten-Nijenhuis bracket of the corresponding
symmetric tensor field on $M$~\cite{norris3}. There is
also a bracket on the direct sum $AHF = \oplus_{p=1} ^\infty AHF^p$ which
gives a graded Poisson algebra, with similar properties.

\subsection{Symplectic vs. $n$-symplectic geometry over $\R^n$}

One interesting point about these vector-valued observables on the frame
bundle is that
their $\R^n$-valued character does not appear when one is looking at their
Poisson algebra. For $T^1$, we can define observables of the form
\[
{\hat \pi_k} = \pi^j _k r_j \qquad {\hat q^i} = q^i r_i ~\mbox {(no sum)}
\qquad {\hat I_k} = r_k
\]
such that
\begin{equation}\label{comm}
\{ {\hat \pi_j}, {\hat q^k} \} = \delta^k _j {\hat I_k} ~\mbox {(no sum)}
\end{equation}
and all other brackets are zero. Aside from the fact there are now $n$ identity
elements, instead of just one, there is no reason to confine our thinking to
$\R^n$ valued functions. Instead, we can consider these as abstract algebraic
objects that satisfy ($\ref{comm}$). One might find it odd that we are only
defining a set of $2n$ observables, since a $n(n+1)$-dimensional manifold
should require
$n(n+1)$ observables to form a complete set. Indeed, we can define such a set
by
\[
{\hat \pi_k} = \pi^j _k r_j \qquad {\hat q^i _j} = q^i r_j
\qquad {\hat I_k} = r_k
\]
where the bracket relation is now
\[
\{ \hpi_k, \hq^i _j \} = \delta^i _k \hI_j
\]
This would give a complete set of observables on $LM$, such that any arbitrary
observable that commutes with these observables is a multiple of the identity.
However, for the functions that satisfy the structure equation of
$n$-symplectic
geometry, it is necessary only to define $2n$ variables, since all the
generalized momenta appear in terms of $\hpi_k$. The functional dimensional of
the $n$-symplectic observables on $LM$ is the same as symplectic observables on
$T^*M$\footnote{One might notice that when $n=1$, then the frame bundle and
the cotangent bundle have the same dimension. If one tries to solve the
structure equation for this case, there is no restriction on the degree of
the generalized momenta, and one recovers symplectic geometry with a $GL(1,\R)
$ scaling action. The $n=1$ case will be excluded when we talk about
quantization, since we already know that the symplectic case cannot be
quantized.}.

From this, it is easy to see that we can map the observables of $LM$
to those of $T^*M$. If we consider $T^*M$ as the associated bundle $LM
\times_{GL(n)} \R^{n*}$, then for $f \in SHF^p$, we can define\footnote{There
are just two $GL(n,\R)$ orbits in $\R^{n*}$, the trivial orbit $\{0\}$ and its
complement in $\R^{n*}$. So, we select $\alpha \neq 0$ to give a map from $LM$
to $T^*M - \{0\}$.}
\begin{equation}
{\tilde f}([u,\alpha]) = \langle f(u), \underbrace{\alpha, \alpha, \cdots,
\alpha}_{p\ \rm{factors}} \rangle
\end{equation}
with the bracket denoting the inner product between elements of $\R^n$ and
$\R^{n*}$, $u = (m, e_j) \in LM$ giving $[u, \alpha] \in T^*M$. First, we see
that
\[
\pi^i _j (m, e_k) \alpha_i = e^i(\partial / \partial x^j) \alpha_i =
p_j(e^i \alpha^i)
\]
so that $\hpi_j \to p_j$. As an example of this map, we see that the
$n$-symplectic observable $f \in ST^2$,
\[
f = f^{ij} (q)\hpi_i \otimes_s \hpi_j
\]
is taken to the symplectic observable
\[
{\tilde f} = f^{ij} (q) p_i p_j
\]
Since we picked
an observable from $ST^p$, the definition is independent of the representative
of the equivalence class. General observables in $SHF^p$ are mapped to
symplectic observables which depend on the choice of representative. For
instance, the observable 
\[
f = f^{ij} (q) \hpi_i \otimes_s \hI_j
\]
in $SHF^2$ is taken to the symplectic observable
\[
{\tilde f} = f^{ij} (q) p_i \alpha_j
\]
This choice can be regarded as a choice of gauge under the structure group.
Since the equivalence classes $[(u,\alpha)]$ are defined by $(u \cdot g,\alpha)
\sim (u, g \cdot \alpha)$, then moving up the fiber of $LM$ changes $\alpha_i$
by a linear transformation to $\alpha_j g^j _i, g \in GL(n,\R)$. This map also
gives the same symplectic observable for many choices of $n$-symplectic
observable. For instance, if we picked a specific choice of gauge for $\alpha$,
and looked at
\[
f = f^{ij} \alpha_j \hpi_i
\]
in $T^1$, we get the same function on the cotangent bundle as when we picked
the element from $SHF^2$ above.

\section{Quantization of observables on $L(\R^n)$}

For the quantization of a symplectic manifold, one usually follows the Dirac
prescription, constructing a map $\cQ$ from observables to operators on a
Hilbert space with the properties
\begin{eqnarray*}
&(1)& [ \cQ(f), \cQ(g) ] = -i\hbar \cQ(\{ f, g \}) \\
&(2)& \cQ(1) = I \\
&(3)& f \to \cQ(f) \mbox{ is a linear map over } \R
\end{eqnarray*}
However, this is often not enough, since it produces a Hilbert space that is
too large. There must be a choice of an additional axiom to the three
above, such as requiring the product of observables
shall be quantized to the Jordan product of operators\footnote{This rule is
to put another set of relations on an overcomplete set of operators. The
problem arises~\cite{ashtekar} when both $F,G$ and their product $FG$ are in
the algebra of observables which one wants to quantize. Then, one must
require some kind of consistency between the operators associated to these
observables. See also~\cite{gotay} for the relation of this Jordan product
rule and the no-go theorem.}, or
\begin{equation}\label{jordan}
\cQ(fg) = {1 \over 2}(\cQ(f) \cQ(g) + \cQ(g) \cQ(f))
\end{equation}
If we look at the abstract properties of the algebra satisfying these rules,
then it can be shown on $\R^n$ to lead to an inconsistency for quartic
observables. This is because, when one seeks to quantize the observable $p^2
q^2$ (using the cotangent bundle of $\R$ as an example), with the Poisson
bracket $\{ p, q\}=1$,
then
\begin{eqnarray*}
\cQ(p^2 q^2) &=& {1 \over 2} (\cQ(p^2) \cQ(q^2) + \cQ(q^2) \cQ(p^2)) \\
&=& {1 \over 2} (\cQ^2(p) \cQ^2(q) + \cQ^2(q) \cQ^2(p)) \\
&=& \cQ^2 (q) \cQ^2(p) - 2i\hbar \cQ(q) \cQ(p) - \hbar^2 I
\end{eqnarray*}
where first we have used the fact that  $(\ref{jordan})$ implies $\cQ(f^2)
= \cQ^2 (f)$ and then used the relation $[\cQ(p), \cQ(q)] = -i\hbar I$. But,
we also have
\begin{eqnarray*}
\cQ(p^2 q^2) &=& {1 \over 4} (\cQ(p) \cQ(q) + \cQ(q) \cQ(p))^2 \\
&=& \cQ^2 (q) \cQ^2 (p) - 2i\hbar \cQ(q) \cQ(p) - {1 \over 4} \hbar^2 I
\end{eqnarray*}
Thus, the two methods of using  $(\ref{jordan})$ give different operators.

Now, the question is whether or not this happens on the frame bundle. We want
to have an algebra of observables on $LM$ to compare with those on the
cotangent bundle, so we use the fact that
\begin{thm}
Let $f \in SHF^p$ and $g \in SHF^q$. Then $f \otimes_s g$ is an observable in
$SHF^{p+q}$.
\end{thm}
We then consider the class of observables that are symmetric tensor products
of $\hpi_a, \hq^b$, and $\hI_c$. So our corresponding function on $LM$ would
be of the form $\hpi_a \otimes_s \hpi_b \otimes_s \hq^c \otimes \hq^d$. 
Modeling our quantization map after the Jordan product in $(\ref{jordan})$,
we use the fact that the symmetric product is defined as
\[
r_a \otimes_s r_b \otimes_s \cdots \otimes_s r_c \equiv r_{(a} \otimes r_b
\otimes \cdots \otimes r_{c)}
\]
to write the quantization map as
\begin{equation}\label{jordan2}
\cQ(f_{i_1} \otimes_s f_{i_2} \otimes_s \cdots \otimes_s f_{i_k}) =
\cQ(f_{(i_1}) \cQ(f_{i_2}) \cdots \cQ(f_{i_k)})
\end{equation}
for the $f_i$ as one of $\hpi_a, \hq^b$, and $\hI_c$
When we quantized $p^2 q^2$ above, the first time used  $(\ref{jordan})$
with $f = p^2$ and $g = q^2$, while the second time used $f=g=pq$. Because the
Jordan product is not associative, this led to a contradiction. This problem
does not arise with the map in $(\ref{jordan2})$ because this map is
associative, so all possible groupings of the $f_i$ give the same operator. So,
for the observable $\hpi_a \otimes_s \hpi_a \otimes_s \hq^a \otimes_s \hq^a$,
with no summing over the index, we have the operator
\[
\cQ(\hpi_a \otimes_s \hpi_a \otimes_s \hq^a \otimes_s \hq^a) = \cQ^2 (\hq^a)
\cQ^2 (\hpi_a) - 2i\hbar \cQ(\hq^a) \cQ(\hpi_a) - {1 \over 3} \hbar^2 \hI_a
\]

The question of whether or not there is an eventual contradiction
along the lines of the symplectic case is unclear at this point. This algebra
has a major difference from the algebra formed by $p_i$ and $q^j$ on $T^*M$.
As alluded to in Section 2.3, for any given observable on $T^*M$,
there are many $n$-symplectic observables that map to it. Each symplectic
observable is "covered" by many $n$-symplectic observables, so the
contradiction in quantization may be lost. When one is finding
the operators for symplectic observables, one gets relations like
\begin{equation}\label{restrict}
\{p^n, pq \} = np^n
\end{equation}
which constrain the operator $\cQ(p^n)$. This does not occur in $n$-symplectic
geometry, since (dropping indices for clarity)
\[
\{ \hpi \otimes_s \hpi, \hpi \otimes_s \hq \} = 2 \hpi \otimes_s \hpi \otimes_s
\hI
\]
Because of this, the operators of $n$-symplectic
observables can be more general than their symplectic counterparts. For
instance, for the observable $\hpi \otimes_s \hpi$, one can add a constant
matrix to the operator obtained by the map above:
\[
\cQ(\hpi \otimes_s \hpi) = \cQ^2 (\hpi) + \Lambda
\]
Additions like this are not permitted in the symplectic quantization because
of restrictions such as  $(\ref{restrict})$. It might be that these
additions are precisely what is needed for a full quantization. In the
proof of the Groenewald-van Hove theorem, the contradiction arises because
the operators for $\{ p^3, q^3\}$ and $3\{p^2 q, pq^2 \}$, which should both
be the operator for $9p^2 q^2$, instead differ by a multiple of the identity.
The $n$-symplectic case gives the same contradiction, but {\it only} if we
keep to the map in  $(\ref{jordan2})$. Adding constant matrices to the
operators
for $\hpi \otimes_s \hpi \otimes_s \hq$ and $\hpi \otimes_s \hq \otimes_s \hq$
whose commutator is the required multiple of the identity would solve the
problem.
We hope to examine the properties of the algebra of $n$-symplectic observables
in future work to see the viability of a full quantization.

So far we have looked at the symmetric observables, which have their
counterparts on $T^*M$. However, there are also the anti-symmetric observables,
which have no analogue in symplectic geometry. We take the wedge products
of the observables $\hpi_a, \hq^b,$ and $\hI_c$, using the anti-symmetric
version of Proposition 1, and use the quantization map
\[
\cQ(f_{i_1} \wedge f_{i_2} \wedge \cdots \wedge f_{i_k}) =
\cQ(f_{[i_1}) \cQ(f_{i_2}) \cdots \cQ(f_{i_k]})
\]
where the brackets represent anti-symmetrization of the indices.
The operators for the quadratic observables are
\[
\cQ(\hpi_a \wedge \hpi_b) = \cQ(\hq^a \wedge \hq^b) = 0 \qquad
\cQ(\hpi_c \wedge \hq^d) = {1 \over 2}\delta^d _c \cQ(\hI_d)
\]
and
\[
\cQ(\hpi_a \wedge \hI_b) = \cQ(\hq^a \wedge \hI_b) = 0
\]
Because the operators for the quadratic observables are either zero, or the
identity operator, the cubic operators are all zero.

\section{The bundle $\cP \simeq \R^2 \times S^1$}

\subsection{$\R^2$-valued observables}

As our first example of a 2-symplectic manifold,
we shall look at the trivial U(1) principal bundle over $\R^2$,
which is isomorphic to $\R^2 \times S^1$. This 2-symplectic manifold can be
thought of as a subbundle of the frame bundle $L(\R^2) \to \R^2$, and we use
this example to illustrate in detail how the formalism might work for a
general principal bundle.  As seen above, the
2-symplectic form on this bundle over $\R^2$ is $\omega = d\theta$, where
\[
\theta = \left(
\begin{array}{c}
\pi^1 _1 dq^1 + \pi^1 _2 dq^2 \\
\pi^2 _1 dq^1 + \pi^2 _2 dq^2
\end{array} \right)
\]
The frame bundle is trivial, so we can reduce the structure group from $GL(2,
\R)$ to any subgroup, in particular to $SO(2) \simeq U(1)$. We denote this
subbundle as $\cP \simeq \R^2 \times S^1$, and write the
generalized momentum in terms of a single coordinate $\phi$, so that
\begin{equation}\label{rotation}
(\pi^i _j) =
\left(
\begin{array}{cc}
\pi^1 _1 & \pi^1 _2 \\
\pi^2 _1 & \pi^2 _2
\end{array}
\right) =
\left(
\begin{array}{cc}
\cos \phi & \sin \phi \\
-\sin \phi & \cos \phi
\end{array}
\right)
\end{equation}
Then, the 2-symplectic potential on the full frame bundle reduces to one on
$\cP$, given by
\begin{equation}\label{oldform}
\theta = \left(
\begin{array}{c}
\theta^1 \\
\theta^2
\end{array}
\right) = \left(
\begin{array}{c}
\cos \phi dq^1 + \sin \phi dq^2 \\
 -\sin \phi dq^1 + \cos \phi dq^2
\end{array}
\right)
\end{equation}
Note that each component is a contact form on $\cP$~\cite{blair}.
We look at how this 2-symplectic form transforms
under the group action. For $SO(2)$, this is given by the translation $\phi
\to \phi + \alpha$; if we use the transformation on the trigonometric
functions, then we see that our 2-symplectic form is also tensorial, and obeys
$R^*_\alpha \theta = R(\alpha) \cdot \theta$, for a rotation $R(\alpha)$ of the
form given in  $(\ref{rotation})$. The subbundle $\cP$ is the bundle of
oriented orthonormal frames on $\R^2$, and $\phi$ gives the angle in $\R^2$
corresponding to the frame.

If we solve the structure equation for vector-valued observables, $f = (f^1,
f^2)$, with $p=1$ and $\omega = d\theta$, then
\[
df^i = - \ X_f \hook \omega^i
\]
where
\[
X_f = X^\phi \pdb{\phi} + X^1 \pdb{q^1} + X^2 \pdb{q^2}
\]
This gives the six equations
\begin{eqnarray}
\sin \phi X^\phi &=& \pd{f^1}{q^1} = \pd{f^2}{q^2}\\
\cos \phi X^\phi &=& - \pd{f^1}{q^2} = \pd{f^2}{q^1}\\
\pd{f^1}{\phi} &=& \cos \phi X^2 - \sin \phi X^1\\
\pd{f^2}{\phi} &=& - \cos \phi X^1 - \sin \phi X^2
\end{eqnarray}
Solving for the Hamiltonian vector field gives
\begin{eqnarray}
X^1 &=& - \sin \phi \pd{f^1}{\phi} - \cos \phi \pd{f^2}{\phi}\\
X^2 &=& \cos \phi \pd{f^1}{\phi} - \sin \phi \pd{f^2}{\phi}
\end{eqnarray}
Thus, only the first four equations place constraints on $(f^1, f^2)$.
Notice that the equations for $X^\phi$ in terms of the observable imply
that the components of $f$ satisfy the Cauchy-Riemann equations. So, if we
introduce the complex variable
$w = q^1 + iq^2$, then $F = f^1 + if^2$ is a function of $w, \phi$
which is holomorphic in $w$. The remaining two equations for
$X^\phi$ give a condition on $F$, namely that
\begin{equation}\label{f-cond}
e^{i\phi} \pd{F}{w} + e^{-i\phi} \pd{\bar F}{\wbar} = 0
\end{equation}
We can also think of this condition as
\[
{\rm Re} \left( e^{i\phi} \pd{F}{w} \right) = 0
\]
Thus, since $F$ is holomorphic in $w$, $(\ref{f-cond})$ implies that $F$ is
linear in $w$, and we can see that the general solution is given by
\[
F = i A(\phi) w e^{-i\phi} + B(\phi) + i C(\phi)
\]
where $A, B,$ and $C$ are all real functions of $\phi$.
By splitting $F$ into its real and imaginary parts, we can see that
it is associated to the observable
\begin{equation}\label{gen1}
f = A(\phi)
\left(
\begin{array}{c}
i w e^{-i\phi} - i\wbar e^{i\phi} \\
w e^{-i\phi} + \wbar e^{i\phi}
\end{array}
\right) + \left(
\begin{array}{c}
B(\phi) \\
C(\phi)
\end{array}
\right)
\end{equation}
Because we have taken the frame bundle of $\R^2$ and reduced the structure
group, we have constrained our observables so that they are affine linear in
$w$. This is in constrast to the case on $LM$, when the vector-valued
observables are affine linear in the generalized momenta, so from that,
one might have expected a restriction on $\phi$.

To find a complete set of observables, we solve for the three basis vectors
of the tangent space. We find three (real) vector fields
\[
X_x = \pdb{q^1} \qquad
X_y = \pdb{q^2}
\]
and
\[
X_p = q^1 \pdb{q^2} - q^2 \pdb{q^1} + \pdb{\phi}
\]
spanning the tangent space of $\cP$, which are the Hamiltonian vector fields
of the three observables
\[
x = \biggl(
\begin{array}{c}
\cos \phi \\
 - \sin \phi
\end{array}
\biggr)
\qquad
y = \biggl(
\begin{array}{c}
\sin \phi \\
\cos \phi
\end{array}
\biggr)
\] and
\[
p = \biggl(
\begin{array}{c}
q^1 \sin \phi - q^2 \cos \phi \\
q^1 \cos \phi + q^2 \sin \phi
\end{array}
\biggr)
\]
These are associated to the complex functions ${\tilde x} = e^{-i\phi}$,
${\tilde y} = ie^{-i\phi}$ and ${\tilde p} = iwe^{-i\phi}$, respectively.
The first two observables generate translations along the base manifold $\R^2$,
while the third generates a rotation of the coordinates of $\R^2$, along with
a translation
in the $U(1)$ degree of freedom. If we define the Poisson bracket by the
formula $\{ f, g \} = X_f (g)$, the three observables obey the following
commutation relations:
\[
\{x,y\} = 0 \qquad \{x,p\} = y \qquad \{y,p\} = -x
\]
These observables form the Lie algebra of the Euclidean group $E(2)$. For two
observables $f,g$, the bracket relation is given in terms of their associated
complex functions ${\tilde f}, {\tilde g}$ as
\[
\{ {\tilde f}, {\tilde g} \} = -i e^{i\phi} \pd{\tilde f}{w}
\pd{\tilde g}{\phi} + i e^{i\phi} \pd{\tilde f}{\phi}
\pd{\tilde g}{w}
\]

When plugging in observables of the form $(\ref{gen1})$ into the bracket
relation, all the important information is carried by the derivatives of the
functions with respect to $\phi$. It seems natural to quantize so that the
linear function quantizes to a derivative operator,
\[
\cQ
\left[ A(\phi)
\left(
\begin{array}{c}
iw e^{-i\phi} - i \wbar e^{i\phi} \\
w e^{-i\phi} + \wbar e^{i\phi}
\end{array}
\right) \right]
=
-i\hbar A(\phi) {\partial \over \partial \phi} \cdot
\left(
\begin{array}{cc}
1& 0 \\
0 & 1
\end{array}
\right)
\]
For the other types of observables, we have simply a multiplication operator:
\begin{eqnarray*}
\cQ \left[ \left(
\begin{array}{c}
B(\phi) \\
0
\end{array}
\right) \right]
&=&
\left(
\begin{array}{cc}
B(\phi) & 0 \\
0 & 0
\end{array}
\right) \\
\cQ \left[ \left(
\begin{array}{c}
0 \\
C(\phi)
\end{array}
\right) \right]
&=&
\left(
\begin{array}{cc}
0 & 0 \\
0 & C(\phi)
\end{array}
\right)
\end{eqnarray*}
A nice feature of this quantization is that all the operators are diagonal.
Note that by replacing the angular variable $\phi$ with the holomorphic
$w$, we can also obtain a quantization map to $\C^2$-valued holomorphic
functions on $\R^2$.

\subsection{Relation to $\C^2$}

Why are we getting this correspondence with complex functions? The reason is
that the 2-symplectic structure of $\cP$ is coming from the symplectic
structure on $\C^2$. To see this, we first note that $\C^2$ is a hyperkahler
manifold, and hence has three complex structures $\I,\J,\K,$ satisfying the
identities of the quaternions\footnote{See~\cite{hyper} for some discussion
of hyperkahler manifolds.},
\[
{\I}^2 = {\J}^2 = {\K}^2 = -1 \qquad {\I\J\K} = -1
\]
Each complex structure has an associated symplectic form, written
in coordinates that are holomorphic with respect to $\I$ as
\begin{eqnarray}\label{forms}
\omega^1 &=& dz \wedge d\zbar + dw \wedge d\wbar \nonumber \\
\omega^2 &=& dz \wedge dw + d\zbar \wedge d\wbar \\
\omega^3 &=& -i dz \wedge dw + i d\zbar \wedge d\wbar \nonumber
\end{eqnarray}
Note that these are the Kaehler forms corresponding to $\I, \J$ and $\K$.
Then we can define holomorphic and anti-holomorphic symplectic forms
$\omega^{\pm}$ for $\I$ by
\[
\omega^{\pm} = \omega^2 \pm i\omega^3
\]
To recoup the 2-symplectic structure on $\cP$, we embed the
manifold in $\C^2$, by the map
\[
(q^1, q^2, \phi) \to (q^1 + iq^2, e^{i\phi})
\]
Then, it is easy to check that the real part of $\omega^+$ is the same as the
differential of $\theta^1$
in $(\ref{oldform})$, and the imaginary part as $d\theta^2$. The
fact that $\omega^+$ is degenerate on $\cP$  is overcome by separating its
real and imaginary parts into a column vector to obtain a 2-symplectic form.

We can also examine the 2-symplectic structures on $\cP$ that are
derived from the $\J$- and $\K$-holomorphic symplectic forms on $\C^2$. It
turns out we find no non-trivial examples of vector observables. As an
example, we look at the 2-symplectic
form given by the $\J$ complex structure on $\C^2$, which is
\[
\omega_{\J} = \left(
\begin{array}{c}
e^{i\phi} d\phi \wedge dw + e^{-i\phi} d\phi \wedge d\wbar \\
dw \wedge d\wbar
\end{array}
\right)
\]
This 2-symplectic form is different from the one previously considered since,
under the group action, it does not transform tensorially. In fact, since the
second component is a volume form on $\R^2$, it is invariant under both the
$U(1)$ of the fiber, and $SO(2)$ rotations of the base manifold.

Setting up the structure equation for vector observables, we find
\[
\pd{f^1}{w} e^{-i\phi} = \pd{f^1}{\wbar} e^{i\phi} \qquad
\pd{f^1}{\phi} = e^{-i\phi} \pd{f^2}{w} - e^{i\phi} \pd{f^2}{\wbar}
\]
as well as $\pd{f^2}{\phi} = 0$. So, if we look for solutions of the form
\[
f = \left(
\begin{array}{c}
\sum _{m} A_m (w,\wbar) e^{im\phi} \\
B(w,\wbar)
\end{array}
\right)
\]
where we have the reality conditions that $B(w,\wbar)$ is real, and $A_{-m} =
A^* _m$. Then the equations on $f$ give that
\[
{\partial A_{m+1} \over \partial w} = {\partial A_{m-1} \over \partial \wbar}
\]
and
\[
i A_{-1} =-{\partial B \over \partial w} \qquad
i A_1 = -{\partial B \over \partial \wbar}
\]
Since we are restricting to only real observables, there are no solutions to
the last two equations beyond a constant function $B$. The fact that the
2-symplectic
form is neither invariant nor tensorial under the $U(1)$ group action has
precluded any kind of non-trivial vector observables.

Finally, we use these 2-symplectic forms on $\C^2$ itself. With no restriction
on the coordinate $z$ now, we have
\begin{equation}\label{2symp}
\omega = \left(
\begin{array}{c}
\omega^2 \\ \omega^3 
\end{array} \right) =
\left(
\begin{array}{c}
dz \wedge dw + d\zbar \wedge d\wbar \\
i d\zbar \wedge d\wbar - i dz \wedge dw
\end{array} \right)
\end{equation}
Again, we have divided the holomorphic symplectic form $\omega^+$ into its
real and imaginary parts. This 2-symplectic form is tensorial under the $U(1)$
transformation, in a manner similar to that of the 2-symplectic form on $\cP$.

Now, we examine which functions can be observables, skipping over details that
are similar to the previous cases. We will show that the $\R^2$ valued
observables of the 2-symplectic geometry on $\C^2$, using the 2-form 
$(\ref{2symp})$ is equivalent to the $\I$-holomorphic symplectic geometry
$(\C^2, \omega^+)$. If we work out the
equations for the global Hamiltonian vector fields from the structure equation
for $p=1$, we get the following equations for $f = (f^1, f^2)$,
\[
{\partial f^1 \over \partial z} = i {\partial f^2 \over \partial z} \qquad
{\partial f^1 \over \partial \zbar} =- i {\partial f^2 \over \partial \zbar}
\]
and similarly for $w$ and $\wbar$,
\[
{\partial f^1 \over \partial w} = i {\partial f^2 \over \partial w} \qquad
{\partial f^1 \over \partial \wbar} =- i {\partial f^2 \over \partial \wbar}
\]
From these equations for the observables, we see that the general form is
\[
f = \left(
\begin{array}{c}
C(z,w) + {\bar C} (\zbar,\wbar) \\
-i C(z,w) + i {\bar C} (\zbar,\wbar)
\end{array}
\right)
\]
or just the real and imaginary parts of a function on $\C^2$ that is
holomorphic with respect to the complex structure $\I$.
 Because of this association, we can
either look at the Poisson brackets of the observables themselves, or their
associated holomorphic functions, where ${\tilde f} = f^1 + if^2$. The latter
is simply the Poisson bracket on $\C^2$, given by
\[
\{ {\tilde f}, {\tilde g} \} = {\partial {\tilde f} \over \partial z}{\partial
{\tilde g} \over \partial w} - {\partial {\tilde f} \over \partial w}{\partial
{\tilde g} \over \partial z}
\]
thus regaining the symplectic structure on $\C^2$.

\section{The 3-sphere}

Our last example is the 3-sphere, which can be thought of as a manifold
embedded in $\C^2$, inducing a 2-symplectic geometry on it. One way to look at
the vector observables on $S^3$ is to
use only those from $\C^2$ whose Hamiltonian vector fields are always tangent
to the 3-sphere. We use the results from above to get the vector fields for
the holomorphic function ${\tilde f} (z,w)$ on $\C^2$ associated to an
observable $f$ as
\[
2X_{\tilde f} = - \pd{\tilde f}{w} \pdb{z} + \pd{\tilde f}{z} \pdb{w} - 
\pd{\bar {\tilde f}}{\wbar} \pdb{\zbar} + \pd{\bar {\tilde f}}{\zbar}
\pdb{\wbar}
\]
 If we use the dot product in $\C^2$ to see which of these
vector fields are perpendicular to the radial vector, then we get the
condition that
\begin{equation}\label{radial2}
\wbar \pd{\tilde f}{z} - \zbar \pd{\tilde f}{w} + w \pd{\bar {\tilde f}}
{\zbar} - z \pd{\bar {\tilde f}}{\wbar} = 0
\end{equation}
The general solution to this constraint equation is ${\tilde f} = C_i x^i$,
where $C_i,\ i=1,2,3,$ are real constants, and
\[
x^1 = -izw \qquad x^2 = {i(z^2 - w^2) \over 2} \qquad x^3 = {z^2 + w^2 \over
2}
\]
We can show this by using a polynomial in $z$ and $w$ for ${\tilde f}$, then
using the condition $(\ref{radial2})$ to limit the coefficients to the above.
If we look at the bracket relations of these variables, we find that
$\{ x^i, x^j \} = \epsilon^{ijk} x^k$, so that these observables give the Lie
algebra of $SU(2)$.

We have naturally gotten observables of the three spin variables $x^i$ on
the 3-sphere. Note that $S^3 \simeq SU(2)$, so that we can also think about
constructing a
vector-valued 2-form on the 3-sphere that is tensorial under $SU(2)$, which
seems a more logical place to obtain these spin variables. Since
$SU(2)$ is a Lie group, the tangent vectors at the identity form a Lie
algebra, with the bracket taking two vectors and giving a third. This can
be thought of as a ${\frak s \frak u}(2)$-valued 2-form on the tangent space of
$SU(2)$, using the group action to map the vectors at any point to the
identity. If we choose a basis of the Lie algebra $\{X_j \}$ and a dual
basis $\{ \theta^k \}$, then this 2-form can be written as
\[
\omega^i = {1 \over 2} c^i _{jk} \theta^j \wedge \theta^k = [\theta, \theta]^i
\]
where $[\ ,\ ]$ denotes the Lie bracket of the group.
However, by the Maurer-Cartan equation, this gives $\omega^i = -d\theta^i$.
Since the Lie algebra of the 3-sphere is isomorphic to $\R^3$, we can also
use the interpretation that this is a $\R^3$-valued 2-form on $S^3$. So,
if we write down a structure equation for vector-valued observables, we have
\begin{equation}\label{s3:struct}
df^i = -X \hook \omega^i = X \hook d\theta^i
\end{equation}

Now we need to pick the basis of vector fields and their dual basis. We use
the fact that the 3-sphere is unit sphere in the space
of quaternions; if we start with the radial vector in $\C^2 = \H$,
\[
r = q^1 {\partial \over \partial q^1}+q^2 {\partial \over \partial q^2}
+ q^3 {\partial \over \partial q^3} + q^4 {\partial \over \partial q^4}
\]
we can associate this with the quaternion ${\tilde q} = q^1 + \i q^2 + \j q^3
+ \k q^4$.
The three complex structures $\I, \J,$ and $\K$ on $\C^2$ are related to the
vectors associated with $\i {\tilde q}, \j {\tilde q},$ and $\k {\tilde q}$,
respectively, given by
\begin{eqnarray*}
v^1  &=& {\I}(r) = - q^2{\partial \over \partial q^1} + q^1{\partial \over
\partial q^2}
- q^4 {\partial \over \partial q^3} + q^3 {\partial \over \partial q^4}
\\
v^2 &=& {\J}(r) = - q^3{\partial \over \partial q^1} + q^4 {\partial \over
\partial q^2}
+ q^1 {\partial \over \partial q^3} - q^2 {\partial \over \partial q^4}
\\
v^3 &=& {\K}(r) = - q^4{\partial \over \partial q^1} - q^3 {\partial \over
\partial q^2}
+ q^2 {\partial \over \partial q^3} + q^1 {\partial \over \partial q^4}
\end{eqnarray*}
These three vectors together form the Lie algebra $\frak s \frak u(2)$, giving
us a global basis of vector fields. The differentials of their dual basis
$\theta^i$ are simply the three symplectic 2-forms from $\C^2$:
\begin{eqnarray*}
d\theta^1 &=& dq^1 \wedge dq^2 + dq^3 \wedge dq^4 \\
d\theta^2 &=& dq^1 \wedge dq^3 - dq^2 \wedge dq^4 \\
d\theta^3 &=& dq^1 \wedge dq^4 + dq^2 \wedge dq^3
\end{eqnarray*}
(these are proportional to the symplectic forms in  $(\ref{forms})$). Together
we use these three 2-forms for our 3-symplectic 2-form, $\omega^i =
-d\theta^i$. Notice that these $\theta^i$ are contact structures on the
3-sphere.

We can see that $v^i \hook d\theta^j = \epsilon^{ijk} \theta^k$ and $\theta^k
\neq df^k$. So, none of these vector fields are Hamiltonian on $S^3$. However,
to find the three vector fields $v^i$, we used the left action of $SU(2)$ on
the radial vector. If we use the right action also, we get the vector fields
\begin{eqnarray*}
w^1  &=& = - q^2{\partial \over \partial q^1} + q^1{\partial \over
\partial q^2}
+ q^4 {\partial \over \partial q^3} - q^3 {\partial \over \partial q^4}
\\
w^2 &=& = - q^3{\partial \over \partial q^1} - q^4 {\partial \over
\partial q^2}
+ q^1 {\partial \over \partial q^3} + q^2 {\partial \over \partial q^4}
\\
w^3 &=& = - q^4{\partial \over \partial q^1} + q^3 {\partial \over
\partial q^2}
- q^2 {\partial \over \partial q^3} + q^1 {\partial \over \partial q^4}
\end{eqnarray*}
These three vector fields also form a basis of the Lie algebra, and, if we
solve ($\ref{s3:struct}$), we find these are also Hamiltonian vector fields
for the three variables
\begin{eqnarray*}
w^1 &\Leftrightarrow& y^1 = {1 \over 2} \bigl[(q^1)^2 + (q^2)^2 - (q^3)^2 -
(q^4)^2 \bigr] \\
w^2 &\Leftrightarrow& y^2 = q^2 q^3 - q^1 q^4 \\
w^3 &\Leftrightarrow& y^3 = q^1 q^3 + q^2 q^4
\end{eqnarray*}
The forms of $y^i$ listed are the first component of the vector-valued
observables,
and, since their Hamiltonian vector fields are a basis of the Lie algebra
of $SU(2)$, then we have that $\{y^i, y^j\} = \epsilon^{ijk} y^k$.
Here again we have obtained observables associated with ${\frak s \frak u(2)}$.
But the situation is different than the $2$-symplectic case -- there, it was
the condition that the Hamiltonian vector fields be tangent to the 3-sphere
that restricted the functions to those of the three spin variables. For the
$3$-symplectic case, there is a similar restriction, but there are also
relations between the derivatives, such as
\[
\pd{f^1}{q^2} = \pd{f^2}{q^3} = \pd{f^3}{q^4}
\]
which further restrict the observables to be linear in the three spin
variables. As with $\cP$, there is a restriction to linearity on
the part of the vector-valued observables.

\section{Discussion}

This paper is intended as a first step to look at observables on
$n$-symplectic manifolds, and at $n$-symplectic geometry as a new model for
both classical dynamics and the quantization procedure. As usual, this has
raised a great many questions,
and there are many directions to be explored, some
of which have been touched on. One area is that of constraints on
$n$-symplectic manifolds. On the frame bundle, observables in $SHF^p$ are
limited to be no more than degree $p$ in the generalized momenta $\pi^i _j$.
However, we seem to have a different story on $\cP$ -- the
vector-valued observables are indeed limited, but are linear in the
variables of the base manifold.
Since we can think of $L(\R^2)$ and $\cP$ as the extreme
cases of principal bundles over $\R^2$, with the largest and smallest
structure groups, one wonders how other groups would affect the possible
observables. We have only looked at the vector-valued observables in this
paper, but in future work we hope to examine how constrainted more general
observables are for various principal bundles.

The exact role of the structure group of the principal bundle is also not
clear. For the 3-sphere, we looked at both a $2$-symplectic form, tensorial
under $U(1)$, as well as a $3$-symplectic form, transforming under $SU(2)$.
Is there any kind of relation between the observables in these theories?
Both give a set of spin observables, so there must be a link between the
two. Note that, in constructing the $3$-symplectic form, we took
the radial vector, acting on it with the left action of the group to give
a basis of the Lie algebra. One could do the same with the right action to
obtain another 2-form. How are these $3$-symplectic forms related? It is
also interesting to consider the case of the 7-sphere. Although it is
not a Lie group, it is a parallelizable manifold, so one could go through with
the same kind of analysis, the only difference being one would act on the
radial vector in $\R^8$ with the octonions instead of the quaternions. Because
there is no structure group associated with this theory, the question remains
of how important is the group in the first place. We compare the 7-sphere to
the case when we
looked at a $2$-symplectic form on $\cP$ which was not tensorial
under $U(1)$, obtaining only constant functions. Also compare with the
definitions
of $n$-symplectic structures of Awane, and De Leon, et al., where there is
no group action on the structure~\cite{awane,deleon1,deleon2}.

Another question is the relation of $n$-symplectic geometry to gauge theory.
For instance, we can think of $\cP$ as a Kaluza-Klein type
manifold. Because of this, one would think the gauge freedom is the usual
transformations on $U(1)$. It was mentioned above that the $n$-symplectic
theory seems to be
mixing the group $U(1)$ with the base manifold $\R^2$ by the gauge 
vector fields occurring only on $\R^2$. Also, if we pick appropriate gauge
conditions, what kind of transformations do we obtain? In a related question,
for the 3-sphere, we used the $SU(2)$ action to formulate $3$-symplectic
geometry. So, we can either think of $S^3$ as the spatial slice of a
spacetime, with $SU(2)$ dynamics on it, or else think of our $3$-symplectic
manifold as a $SU(2)$ principal bundle over a point. With this viewpoint,
we can pick a manifold $M$ with a 3-sphere over each point, giving us
a $3$-symplectic form obtained by integrating over the manifold,
\[
\omega (X, Y) = \int_M \ [X(q), Y (q)]\  d^n q
\]
where $X$ and $Y$ are ${\frak s \frak u}(2)$-valued vectors on $M$.
Note the contrast with typical field theories, where the trace is taken of the
two vectors in the Lie algebra. This procedure can be generalized for any
Lie group which has structure constants $c^i _{jk}$ which are non-degenerate
as a matrix, i.e. for any semi-simple Lie group.

Finally, can $n$-symplectic geometry predict the no-go theorems of geometric
quantization of symplectic manifolds? As we have seen, the observables on $LM$
cover those on $T^*M$, so that the quantization operators on the cotangent
bundle are in some sense induced from those on the frame bundle. Since we can
also consider the 3-sphere as a $U(1)$ bundle over $S^2$, then the three
spin variables and their symplectic geometry on the 2-sphere also can be
studied in a $n$-symplectic context. The symplectic quantization fails at the
level of cubic observables on $T^*\R^n$ but at quadratics on $S^2$. Can this
be related to their $n$-symplectic counterparts, by giving a general
result of when the symplectic quantization will fail?

\section{Acknowledgements}

The author would like to thank Ranee Brylinski and Larry Norris for many
helpful discussions.

\end{document}